\documentclass[12pt]{article}
\usepackage{amssymb}
\usepackage{amsfonts}
\usepackage{amsmath}

\input{tcilatex}
\begin{document}

\begin{center}
\textbf{New fixed point theorems for set-valued contractions }

\textbf{in }$b$\textbf{-metric spaces}

\bigskip

Radu MICULESCU and Alexandru MIHAIL

\bigskip
\end{center}

\textbf{Abstract}. {\small In this paper we indicate a way to generalize a
series of fixed point results in the framework of }$b${\small -metric spaces
and we exemplify it by extending Nadler's contraction principle for
set-valued functions (see Multi-valued contraction mappings, Pac. J. Math., 
\textbf{30} (1969), 475-488) and a fixed point theorem for set-valued
quasi-contractions functions due to H. Aydi, M.F. Bota, E. Karapinar and S.
Mitrovi\'{c} (see A fixed point theorem for set-valued quasi-contractions in 
}$b${\small -metric spaces, Fixed Point Theory Appl. 2012, 2012:88).}

\bigskip



\bigskip

\textbf{1.} \textbf{Introduction}

\bigskip

In the last decades one can observe a remarkable amount of interest for the
development of fixed point theory since it has a huge number of applications.

Among the generalizations of the Banach-Caccioppoli-Picard principle - one
of the central results of the above mentioned theory, known also as the
contraction principle - a central role is played by the following two:

- the one due to S.B. Nadler [21] who extended the contraction principle to
set-valued functions and generated in this way many applications in control
theory, convex optimization etc;

- the one due to I. A. Bakhtin [5] and S. Czerwik [13], [14] who, motivated
by the problem of the convergence of measurable functions with respect to
measure, introduced $b$-metric spaces (a generalization of metric spaces)
and proved the contraction principle in this framework.

In the last period many mathematicians obtained fixed point results for
single-valued or set-valued functions, in the setting of $b$-metric spaces
(see, for example, [1], [8], [9], [10], [17], [23], [24], [29], [30] and the
references therein).

In this paper we indicate a way (see Lemma 2.2) to generalize a series of
fixed point results in the framework of $b$-metric spaces and we exemplify
it by extending{\small \ }Nadler's contraction principle for set-valued
functions (see [21]) and a fixed point theorem for set-valued
quasi-contractions functions due to H. Aydi, M.F. Bota, E. Karapinar and S.
Mitrovi\'{c} (see [4]).

\bigskip

\textbf{2.} \textbf{Preliminaries results}

\bigskip

In this section we sum up some basic facts that we are going to use later.

\bigskip

\textbf{Definition 2.1.} \textit{Given a nonempty set }$X$\textit{\ and a
real number }$s\in \lbrack 1,\infty )$\textit{, a function }$d:X\times
X\rightarrow \lbrack 0,\infty )$\textit{\ is called }$b$\textit{-metric if
it satisfies the following properties:}

\textit{i) }$d(x,y)=0$\textit{\ if and only if }$x=y$\textit{;}

\textit{ii) }$d(x,y)=d(y,x)$\textit{\ for all }$x,y\in X$\textit{;}

\textit{iii) }$d(x,y)\leq s(d(x,z)+d((z,y))$\textit{\ for all }$x,y,z\in X$%
\textit{.}

\textit{The pair }$(X,d)$\textit{\ is called }$b$\textit{-metric space.}

\bigskip

\textbf{Remark 2.1}. \textit{As when }$s=1$\textit{, a }$b$\textit{-metric
space is a metric space, we infer that the family of }$b$\textit{-metric
spaces is larger than the one of metric spaces. In other words, every metric
spaces is a }$b$\textit{-metric space. Note that Czerwik proved that the
converse need not be true (see also [4], [12], [18], [22] and [27]), so the
family of }$b$\textit{-metric spaces is effectively larger than the one of
metric spaces.}

\bigskip

\textbf{Definition 2.2.} \textit{A sequence} $(x_{n})_{n\in \mathbb{N}}$ 
\textit{of elements from a }$b$\textit{-metric space }$(X,d)$\textit{\ is
called:}

\textit{- convergent if there exists }$l\in \mathbb{R}$\textit{\ such that }$%
\underset{n\rightarrow \infty }{\lim }d(x_{n},l)=0$;

\textit{- Cauchy if }$\underset{m,n\rightarrow \infty }{\lim }%
d(x_{m},x_{n})=0$\textit{, i.e. for every }$\varepsilon >0$\textit{\ there
exists} $n_{\varepsilon }\in \mathbb{N}$ \textit{such }$d(x_{m},x_{n})<%
\varepsilon $\textit{\ that for all} $m,n\in \mathbb{N}$, $m,n\geq
n_{\varepsilon }$.

\textit{The} $b$\textit{-metric space }$(X,d)$ \textit{is called complete if
every Cauchy sequence of elements from }$(X,d)$\textit{\ is convergent.}

\bigskip

Beside the classical spaces $l^{p}(\mathbb{R})$ and $L^{p}[0,1]$, where $%
p\in (0,1)$, one can find examples of $b$-metric spaces in [4], [6], [10],
[13] and [14].

\bigskip

\textbf{Remark 2.2.} \textit{As in the case of metric spaces, a }$b$\textit{%
-metric space can be endowed with the topology induced by its convergence
and almost all the concepts and results which are valid for metric spaces
can be extended to the framework of }$b$\textit{-metric spaces.}

\textit{T.V. An, L.Q. Tuyen and N.V. Dung [3] proved that every }$b$\textit{%
-metric space is a semi-metrizable space (i.e. there exists a function }$%
d:X\times X\rightarrow \lbrack 0,\infty )$\textit{\ such that: i) }$d(x,y)=0$%
\textit{\ if and only if }$x=y$\textit{; ii) }$d(x,y)=d(y,x)$\textit{\ for
all }$x,y\in X$\textit{; iii) }$x\in \overline{A}$ \textit{if and only if }$%
d(x,A)=\inf \{d(x,y)\mid y\in A\}=0$\textit{\ for every }$x\in X$\textit{\
and every }$A\subseteq X$\textit{). Consequently many properties of }$b$%
\textit{-metric spaces are obvious. In addition they provided a sufficient
condition for a }$b$\textit{-metric space to be metrizable and gave an
example showing that, in the framework of a }$b$\textit{-metric space }$%
(X,d) $\textit{, there exists an open ball (i.e. a set of the form }$\{y\in
X\mid d(x,y)<r\}$\textit{, where }$r>0$\textit{) which is not open.}

\textit{In a metric space }$(X,d)$\textit{, the functions }$d$\textit{\ is
continuous (i.e. }$\underset{n\rightarrow \infty }{\lim }%
d(x_{n},y_{n})=d(x,y)$\textit{\ for all sequences} $(x_{n})_{n\in \mathbb{N}%
} $ \textit{and} $(y_{n})_{n\in \mathbb{N}}$ \textit{of elements from }$X$%
\textit{\ and }$x,y\in X$\textit{\ such that }$\underset{n\rightarrow \infty 
}{\lim }x_{n}=x$\textit{\ and }$\underset{n\rightarrow \infty }{\lim }%
y_{n}=y $). \textit{The fact that this property is not valid for }$b$\textit{%
-metric spaces (as }$\frac{1}{s^{2}}d(x,y)\leq \underset{n\rightarrow \infty 
}{\underline{\lim }}d(x_{n},y_{n})\leq \underset{n\rightarrow \infty }{%
\overline{\lim }}d(x_{n},y_{n})\leq s^{2}d(x,y)$ \textit{and }$\frac{1}{s}%
d(x,y)\leq \underset{n\rightarrow \infty }{\underline{\lim }}d(x_{n},y)\leq 
\underset{n\rightarrow \infty }{\overline{\lim }}d(x_{n},y)\leq sd(x,y)$, 
\textit{see [20], [22] and [25]) is a motivation of our Definition 3.2.}

\textit{\bigskip }

In the sequel, given a $b$-metric space $(X,d)$:

- by $\mathcal{B}(X)$ we denote the set of non-empty bounded closed subsets
of $X$

- for $A,B\in \mathcal{B}(X)$, we define the Hausdorff-Pompeiu distance
between $A$ and $B$ by $h(A,B)=\max \{\underset{a\in A}{\sup }d(a,B),%
\underset{b\in B}{\sup }d(b,A)\}$, where $d(x,C)=\underset{c\in C}{\inf }%
d(x,c)$ for every $x\in X$ and every $C\in \mathcal{B}(X)$

- for $c,d\in \lbrack 0,1]$ and $x,y\in X$, we shall use the following
notation:

\begin{equation*}
N_{c,d}(x,y)=\max \{d(x,y),cd(x,T(x)),cd(y,T(y)),\frac{d}{2}%
(d(x,T(y))+d(y,T(x)))\}\text{.}
\end{equation*}

\bigskip

\textbf{Lemma 2.1}. \textit{For every sequence} $(x_{n})_{n\in 
%TCIMACRO{\U{2115} }%
%BeginExpansion
\mathbb{N}
%EndExpansion
}$ \textit{of elements from a }$b$\textit{-metric space }$(X,d)$\textit{,
the inequality }%
\begin{equation*}
d(x_{0},x_{k})\leq s^{n}\underset{i=0}{\overset{k-1}{\sum }}d(x_{i},x_{i+1})
\end{equation*}%
\textit{is valid for every }$n\in \mathbb{%
%TCIMACRO{\U{2115} }%
%BeginExpansion
\mathbb{N}
%EndExpansion
}$\textit{\ and every }$k\in \{1,2,3,...,2^{n}-1,2^{n}\}$\textit{.}

\textit{Proof}. We are going to use the method of mathematical induction.
Denoting by $P(n)$ the statement: $d(x_{0},x_{k})\leq s^{n}\underset{i=0}{%
\overset{k}{\sum }}d(x_{i},x_{i+1})$ for every $n\in 
%TCIMACRO{\U{2115} }%
%BeginExpansion
\mathbb{N}
%EndExpansion
$ and every $k\in \{1,2,3,...,2^{n}-1,2^{n}\}$, as the statements $P(0)$ and 
$P(1)$ are obvious, it remains to prove that $P(n)\Rightarrow P(n+1)$.

Indeed, the above mentioned implication is true since, on the one hand, for
every $k\in \{1,2,3,...,2^{n}-1,2^{n}\}$, we have 
\begin{equation*}
d(x_{0},x_{k})\overset{P(n)}{\leq }s^{n}\underset{i=0}{\overset{k}{\sum }}%
d(x_{i},x_{i+1})\leq s^{n+1}\underset{i=0}{\overset{k}{\sum }}%
d(x_{i},x_{i+1})\text{.}
\end{equation*}%
On the other hand, for every $k\in \{2^{n}+1,2^{n}+2,...,2^{n+1}-1,2^{n+1}\}$%
, we have 
\begin{equation*}
d(x_{0},x_{k})\leq s(d(x_{0},x_{2^{n}})+d(x_{2^{n}},x_{k}))\overset{P(n)}{%
\leq }
\end{equation*}%
\begin{equation*}
\leq s(s^{n}\underset{i=0}{\overset{2^{n}-1}{\sum }}d(x_{i},x_{i+1})+s^{n}%
\underset{i=2^{n}}{\overset{k-1}{\sum }}d(x_{i},x_{i+1}))=s^{n+1}\underset{%
i=0}{\overset{k-1}{\sum }}d(x_{i},x_{i+1})\text{. }\square
\end{equation*}

\bigskip

\textbf{Lemma 2.2}. \textit{Every sequence} $(x_{n})_{n\in 
%TCIMACRO{\U{2115} }%
%BeginExpansion
\mathbb{N}
%EndExpansion
}$ \textit{of elements from a }$b$\textit{-metric space }$(X,d)$\textit{,
having the property that there exists }$\gamma \in \lbrack 0,1)$ \textit{%
such that }%
\begin{equation*}
d(x_{n+1},x_{n})\leq \gamma d(x_{n},x_{n-1})\text{,}
\end{equation*}%
\textit{\ for every }$n\in \mathbb{%
%TCIMACRO{\U{2115} }%
%BeginExpansion
\mathbb{N}
%EndExpansion
}$\textit{,} \textit{is Cauchy.}

\textit{Proof}. First let us note that%
\begin{equation}
d(x_{n+1},x_{n})\leq \gamma ^{n}d(x_{1},x_{0})\text{,}  \tag{1}
\end{equation}%
for every $n\in \mathbb{%
%TCIMACRO{\U{2115} }%
%BeginExpansion
\mathbb{N}
%EndExpansion
}$\textit{.}

For all $m,k\in \mathbb{%
%TCIMACRO{\U{2115} }%
%BeginExpansion
\mathbb{N}
%EndExpansion
}$, with the notation $p=[\log _{2}k]$, we have%
\begin{equation*}
d(x_{m+1},x_{m+k})\leq sd(x_{m+1},x_{m+2})+sd(x_{m+2},x_{m+k})\leq
\end{equation*}%
\begin{equation*}
\leq
sd(x_{m+1},x_{m+2})+s^{2}d(x_{m+2},x_{m+2^{2}})+s^{2}d(x_{m+2^{2}},x_{m+k})%
\leq
\end{equation*}%
\begin{equation*}
\leq
sd(x_{m+1},x_{m+2})+s^{2}d(x_{m+2},x_{m+2^{2}})+s^{3}d(x_{m+2^{2}},x_{m+2^{3}})+s^{3}d(x_{m+2^{3}},x_{m+k})\leq
\end{equation*}%
\begin{equation*}
...
\end{equation*}%
\begin{equation*}
\leq \overset{p}{\underset{n=1}{\sum }}%
s^{n}d(x_{m+2^{n-1}},x_{m+2^{n}})+s^{p+1}d(x_{m+2^{p}},x_{m+k})\overset{%
\text{Lemma 1}}{\leq }
\end{equation*}%
\begin{equation*}
\leq \overset{p}{\underset{n=1}{\sum }}s^{2n}(\underset{i=m}{\overset{%
m+2^{n-1}-1}{\sum }}d(x_{2^{n-1}+i},x_{2^{n-1}+i+1}))+s^{2(p+1)}(\underset{%
i=m}{\overset{m+k-2^{p}-1}{\sum }}d(x_{2^{p}+i},x_{2^{p}+i+1}))\leq
\end{equation*}%
\begin{equation*}
\leq \overset{p+1}{\underset{n=1}{\sum }}s^{2n}(\underset{i=m}{\overset{%
m+2^{n-1}-1}{\sum }}d(x_{2^{n-1}+i},x_{2^{n-1}+i+1}))\overset{(1)}{\leq }%
d(x_{0},x_{1})\overset{p+1}{\underset{n=1}{\sum }}s^{2n}(\underset{i=0}{%
\overset{2^{n-1}-1}{\sum }}\gamma ^{m+2^{n-1}+i})\leq
\end{equation*}%
\begin{equation*}
\leq \frac{d(x_{0},x_{1})\gamma ^{m}}{1-\gamma }\overset{p+1}{\underset{n=1}{%
\sum }}s^{2n}\gamma ^{2^{n-1}}=\gamma ^{m}\frac{d(x_{0},x_{1})}{1-\gamma }%
\overset{p+1}{\underset{n=1}{\sum }}\gamma ^{2n\log _{\gamma }s+2^{n-1}}%
\text{.}
\end{equation*}%
Let us note that since $\underset{n\rightarrow \infty }{\lim }(2n\log
_{\gamma }s+2^{n-1}-n)=\infty $, there exists $n_{0}\in \mathbb{N}$ such
that $2n\log _{\gamma }s+2^{n-1}-n\geq M$, i.e. $\gamma ^{2n\log _{\gamma
}s+2^{n-1}}\leq \gamma ^{M}\gamma ^{n}$ for each $n\in \mathbb{N}$, $n\geq
n_{0}$, hence the series $\overset{\infty }{\underset{n=1}{\sum }}\gamma
^{2n\log _{\gamma }s+2^{n-1}}$ is convergent and denoting by $S$ its sum, we
come to the conclusion that%
\begin{equation*}
d(x_{m+1},x_{m+k})\leq \gamma ^{m}\frac{d(x_{0},x_{1})S}{1-\gamma }\text{,}
\end{equation*}%
for all $m,k\in \mathbb{%
%TCIMACRO{\U{2115} }%
%BeginExpansion
\mathbb{N}
%EndExpansion
}$. Consequently, as $\underset{n\rightarrow \infty }{\lim }\gamma ^{n}=0$,
we infer that $(x_{n})_{n\in 
%TCIMACRO{\U{2115} }%
%BeginExpansion
\mathbb{N}
%EndExpansion
}$ is Cauchy. $\square $

\bigskip

\textbf{Theorem 2.1}. \textit{Let }$(X,d)$\textit{\ be a }$b$\textit{-metric
space and} $T:X\rightarrow \mathcal{B}(X)$ \textit{having the property that
there exist }$c,d\in \lbrack 0,1]$ \textit{and} $\alpha \in \lbrack 0,1)$ 
\textit{such that:}

\textit{i) }$\alpha ds<1$\textit{;}

\textit{ii) }$h(T(x),T(y))\leq \alpha N_{c,d}(x,y)$\textit{\ for all }$%
x,y\in X$\textit{.}

\textit{Then, for every }$x_{0}\in X$\textit{, there exist }$\gamma \in
\lbrack 0,1)$ \textit{and a sequence }$(x_{n})_{n\in \mathbb{N}}$ \textit{of
elements from }$X$\textit{\ such that:}

\textit{a) }$x_{n+1}\in T(x_{n})$\textit{\ for every }$n\in \mathbb{N}$%
\textit{;}

\textit{b)} $d(x_{n+1},x_{n})\leq \gamma d(x_{n},x_{n-1})$ \textit{for every 
}$n\in \mathbb{N}$\textit{;}

\textit{c)} $(x_{n})_{n\in \mathbb{N}}$ \textit{is Cauchy.}

\textit{Proof}. Let us consider $\beta \in (\alpha ,\min (1,\frac{1}{ds}))$, 
$\gamma =\max \{\beta ,\frac{ds\beta }{2-ds\beta }\}<1$, $x_{0}\in X$ and $%
x_{1}\in T(x_{0})$.

If $x_{1}=x_{0}$, then the sequence $(x_{n})_{n\in \mathbb{N}}$ given by $%
x_{n}=x_{0}$ for every $n\in \mathbb{N}$ satisfies a), b) and c).

Since $d(x_{1},T(x_{1}))\leq d(T(x_{0}),T(x_{1}))\leq h(T(x_{0}),T(x_{1}))%
\overset{ii)}{\leq }\alpha N_{c,d}(x_{0},x_{1})<\beta N_{c,d}(x_{0},x_{1})$,
there exists $x_{2}\in T(x_{1})$ such that $d(x_{1},x_{2})<\beta
N_{c,d}(x_{0},x_{1})$.

If $x_{2}=x_{1}$, then the sequence $(x_{n})_{n\in \mathbb{N}}$ given by $%
x_{n}=x_{1}$ for every $n\in \mathbb{N}$, $n\geq 1$, satisfies a), b) and c).

By repeating this procedure we obtain a sequence $(x_{n})_{n\in \mathbb{N}}$
of elements from $X$\ such that $x_{n+1}\in T(x_{n})$ and $%
0<d(x_{n},x_{n+1})<\beta N_{c,d}(x_{n-1},x_{n})$ for every\textit{\ }$n\in 
\mathbb{N}$, $n\geq 1$.

With the notation $d_{n}=d(x_{n},x_{n+1})$, we have 
\begin{equation*}
0<d_{n}<\beta N_{c,d}(x_{n-1},x_{n})=
\end{equation*}%
\begin{equation*}
=\beta \max \{d_{n-1},cd(x_{n-1},T(x_{n-1})),cd(x_{n},T(x_{n})),\frac{d}{2}%
(d(x_{n-1},T(x_{n}))+d(x_{n},T(x_{n-1})))\}\leq
\end{equation*}%
\begin{equation*}
\leq \beta \max \{d_{n-1},cd_{n},cd_{n-1},\frac{d}{2}d(x_{n-1},x_{n+1})\}%
\leq \beta \max \{d_{n-1},cd_{n},cd_{n-1},\frac{ds}{2}(d_{n-1}+d_{n})\}\leq
\end{equation*}%
\begin{equation*}
\leq \beta \max \{d_{n-1},\frac{ds}{2}(d_{n-1}+d_{n})\}\text{,}
\end{equation*}%
for every\textit{\ }$n\in \mathbb{N}$, where the justification of the last
inequality is the following: if, by reduction ad absurdum, $\max
\{d_{n-1},cd_{n},cd_{n-1},\frac{ds}{2}(d_{n-1}+d_{n})\}=cd_{n}$, then we get
that $0<d_{n}<\beta cd_{n}\leq \beta d_{n}$, so we obtain the contradiction $%
1<\beta $.

Consequently $d_{n}<\beta d_{n-1}$ or $d_{n}<\beta \frac{ds}{2}%
(d_{n-1}+d_{n})$, i.e. $d_{n}<\beta d_{n-1}$ or $d_{n}<\frac{ds\beta }{%
2-ds\beta }d_{n-1}$ for every\textit{\ }$n\in \mathbb{N}$. Thus $d_{n}\leq
\max \{\beta ,\frac{ds\beta }{2-ds\beta }\}d_{n-1}$, i.e. $%
d(x_{n+1},x_{n})\leq \gamma d(x_{n},x_{n-1})$ for every\textit{\ }$n\in 
\mathbb{N}$.

Hence the sequence $(x_{n})_{n\in \mathbb{N}}$ satisfies a) and b). From
Lemma 2.2 we deduce that it also satisfies c). $\square $

\bigskip

\textbf{3.} \textbf{Main results}

\bigskip

In this section, making use of Theorem 2.1, we present three fixed point
theorems for set-valued functions.

\bigskip

\textbf{Definition 3.1}. \textit{A function} $T:X\rightarrow \mathcal{B}(X)$%
\textit{, where }$(X,d)$\textit{\ is a }$b$\textit{-metric space, is called
continuous if for all sequences }$(x_{n})_{n\in \mathbb{N}}$ \textit{and }$%
(y_{n})_{n\in \mathbb{N}}$ \textit{of elements from} $X$ \textit{and }$%
x,y\in X$\textit{\ such that }$\underset{n\rightarrow \infty }{\lim }x_{n}=x$%
\textit{,} $\underset{n\rightarrow \infty }{\lim }y_{n}=y$ \textit{and} $%
y_{n}\in T(x_{n})$ \textit{for every} $n\in \mathbb{N}$\textit{, we have }$%
y\in T(x)$\textit{.}

\bigskip

\textbf{Theorem 3.1}. \textit{A function} $T:X\rightarrow \mathcal{B}(X)$%
\textit{, where }$(X,d)$\textit{\ is a complete }$b$\textit{-metric space,
has a fixed point, provided that it satisfies the following three conditions:%
}

\textit{i) }$T$ \textit{is continuous;}

\textit{ii) there exist }$c,d\in \lbrack 0,1]$ \textit{and }$\alpha \in
\lbrack 0,1)$\textit{\ such that} $h(T(x),T(y))\leq \alpha N_{c,d}(x,y)$ 
\textit{for all }$x,y\in X$\textit{;}

\textit{iii) }$\alpha ds<1$\textit{.}

\textit{Proof}. Taking into account ii) and iii), by virtue of Theorem 2.1,
there exists a Cauchy sequence $(x_{n})_{n\in \mathbb{N}}$ of elements of $X$%
\ such that 
\begin{equation}
x_{n+1}\in T(x_{n})\text{,}  \tag{1}
\end{equation}%
\ for every $n\in \mathbb{N}$.

As the $b$-metric space $(X,d)$ is complete, there exists $u\in X$ such that 
$\underset{n\rightarrow \infty }{\lim }x_{n}=u$ (so $\underset{n\rightarrow
\infty }{\lim }x_{n+1}=u$). We combine i) with $1)$ to see that $u\in T(u)$,
i.e. $u$ is a fixed point of $T$. $\square $

\bigskip

\textbf{Definition 3.2}. \textit{Given a }$b$\textit{-metric space }$(X,d)$%
\textit{, the }$b$-\textit{metric }$d$ \textit{is called} $\ast $-\textit{%
continuous if for every} $A\in \mathcal{B}(X)$\textit{, every} $x\in X$ 
\textit{and every sequence }$(x_{n})_{n\in \mathbb{N}}$ \textit{of elements
from} $X$ \textit{\ such that }$\underset{n\rightarrow \infty }{\lim }%
x_{n}=x $\textit{, we have }$\underset{n\rightarrow \infty }{\lim }%
d(x_{n},A)=d(x,A)$\textit{.}

\bigskip

Our notion of $\ast $-continuity of $d$ is stronger than the continuity of $%
d $ in the first variable.

\bigskip

\textbf{Theorem 3.2}. \textit{A function} $T:X\rightarrow \mathcal{B}(X)$%
\textit{, where }$(X,d)$\textit{\ is a complete }$b$\textit{-metric space,
has a fixed point, provided that it satisfies the following three conditions:%
}

\textit{i) }$d$ \textit{is }$\ast $-\textit{continuous;}

\textit{ii) there exist }$c,d\in \lbrack 0,1]$ \textit{and }$\alpha \in
\lbrack 0,1)$\textit{\ such that} $h(T(x),T(y))\leq \alpha N_{c,d}(x,y)$ 
\textit{for all }$x,y\in X$\textit{;}

\textit{iii) }$\alpha ds<1$\textit{.}

\textit{Proof}. Based on ii) and iii), according to Theorem 2.1, there
exists a Cauchy sequence $(x_{n})_{n\in \mathbb{N}}$ of elements of $X$\
such that 
\begin{equation}
x_{n+1}\in T(x_{n})\text{,}  \tag{1}
\end{equation}%
\ for every $n\in \mathbb{N}$.

As the $b$-metric space $(X,d)$ is complete, there exists $u\in X$ such that 
$\underset{n\rightarrow \infty }{\lim }x_{n}=u$.

Then we have%
\begin{equation*}
d(x_{n+1},T(u))\overset{(1)}{\leq }d(T(x_{n}),T(u))\leq h(T(x_{n}),T(u))%
\overset{i)}{\leq }
\end{equation*}%
\begin{equation*}
\leq \alpha N_{c,d}(x_{n},u)=
\end{equation*}%
\begin{equation*}
=\alpha \max \{d(x_{n},u),cd(x_{n},T(x_{n})),cd(u,T(u)),\frac{d}{2}%
(d(x_{n},T(u))+d(u,T(x_{n})))\}\overset{(1)}{\leq }
\end{equation*}%
\begin{equation}
\leq \alpha \max \{d(x_{n},u),cd(x_{n},x_{n+1}),cd(u,T(u)),\frac{d}{2}%
(s(d(x_{n},u)+d(u,T(u)))+d(u,x_{n+1}))\}\text{,}  \tag{2}
\end{equation}%
for every $n\in \mathbb{N}$.

Since $\underset{n\rightarrow \infty }{\lim }d(x_{n},u)=\underset{%
n\rightarrow \infty }{\lim }d(x_{n},x_{n+1})=0$ and\ $\underset{n\rightarrow
\infty }{\lim }d(x_{n+1},T(u))=d(u,T(u))$ (as $d$ is $\ast $-continuous and $%
\underset{n\rightarrow \infty }{\lim }x_{n+1}=u$), upon passing to limit, as 
$n\rightarrow \infty $, in $2)$, we get%
\begin{equation}
d(u,T(u)\leq \alpha \max \{0,cd(u,T(u)),\frac{ds}{2}d(u,T(u))\}\leq \max
\{\alpha c,\frac{\alpha ds}{2}\}d(u,T(u))\text{.}  \tag{3}
\end{equation}

As $\max \{\alpha c,\frac{\alpha ds}{2}\}<1$ (see iii)), from $3)$, we
conclude that $d(u,T(u))=0$, i.e. $u\in T(u)$. Hence $T$ has a fixed point. $%
\square $

\bigskip

\textbf{Theorem 3.3}. \textit{A function} $T:X\rightarrow \mathcal{B}(X)$%
\textit{, where }$(X,d)$\textit{\ is a complete }$b$\textit{-metric space,
has a fixed point, provided that it satisfies the following two conditions:}

\textit{i) there exist }$c,d\in \lbrack 0,1]$ \textit{and }$\alpha \in
\lbrack 0,1)$\textit{\ such that} $h(T(x),T(y))\leq \alpha N_{c,d}(x,y)$ 
\textit{for all }$x,y\in X$\textit{;}

\textit{ii) }$\max \{\alpha cs,\alpha ds\}<1$\textit{.}

\textit{Proof}. Making use of i) and ii), according to Theorem 2.1, there
exists a Cauchy sequence $(x_{n})_{n\in \mathbb{N}}$ of elements from $X$\
such that $x_{n+1}\in T(x_{n})$,\ for every $n\in \mathbb{N}$. As the $b$%
-metric space $(X,d)$ is complete, there exists $u\in X$ such that $\underset%
{n\rightarrow \infty }{\lim }x_{n}=u$.

First let us note that, as we have seen in $2)$ from the proof of Theorem
3.2, we have%
\begin{equation*}
d(x_{n+1},T(u))\leq
\end{equation*}%
\begin{equation*}
\leq \alpha \max \{d(x_{n},u),cd(x_{n},T(x_{n})),cd(u,T(u)),\frac{d}{2}%
(d(x_{n},T(u))+d(u,T(x_{n})))\}\leq
\end{equation*}%
\begin{equation*}
\leq \alpha \max \{d(x_{n},u),cd(x_{n},T(x_{n})),cd(u,T(u)),\frac{d}{2}%
(d(x_{n},T(u))+d(u,x_{n+1}))\}\leq
\end{equation*}%
\begin{equation}
\leq \alpha \max \{d(x_{n},u),cd(x_{n},x_{n+1}),cd(u,T(u)),\frac{d}{2}%
(s(d(x_{n},u)+d(u,T(u)))+d(u,x_{n+1})))\}\text{,}  \tag{1}
\end{equation}%
for every $n\in \mathbb{N}$.

We divide the discussion into two cases:

\qquad A. $d(u,T(u))\leq \overline{\underset{n\rightarrow \infty }{\lim }}%
d(x_{n},T(u))$;

and

\qquad B. $d(u,T(u))>\overline{\underset{n\rightarrow \infty }{\lim }}%
d(x_{n},T(u))$.

In case A, there exists a subsequence $(x_{n_{k}})_{k\in \mathbb{N}}$ of $%
(x_{n})_{n\in \mathbb{N}}$ having the property that\ $\underset{k\rightarrow
\infty }{\lim }d(x_{n_{k}+1},T(u))\geq d(u,T(u))$, so for every $\varepsilon
>0$ there exists $k_{\varepsilon }\in \mathbb{N}$ such that 
\begin{equation*}
d(u,T(u))-\varepsilon \leq d(x_{n_{k}+1},T(u))\overset{(1)}{\leq }
\end{equation*}%
\begin{equation*}
\leq \alpha \max \{d(x_{n_{k}},u),cd(x_{n_{k}},x_{n_{k}+1}),cd(u,T(u)),\frac{%
d}{2}(s(d(x_{n_{k}},u)+d(u,T(u)))+d(u,x_{n_{k}+1})))\}\text{,}
\end{equation*}%
for every $k\in \mathbb{N}$, $k\geq k_{\varepsilon }$. By passing to limit
as $k\rightarrow \infty $ in the above inequality, we get that%
\begin{equation*}
d(u,T(u))-\varepsilon \leq \alpha \max \{0,0,cd(u,T(u)),\frac{sd}{2}%
d(u,T(u))\}=d(u,T(u))\max \{\alpha c,\frac{\alpha sd}{2}\}\text{,}
\end{equation*}%
for every $\varepsilon >0$, so 
\begin{equation*}
d(u,T(u))\leq d(u,T(u))\max \{\alpha c,\frac{\alpha sd}{2}\}\text{.}
\end{equation*}%
Since $\max \{\alpha c,\frac{\alpha sd}{2}\}<1$, from the above inequality,
we conclude that\linebreak\ $d(u,T(u))=0$, i.e. $u\in T(u)$, so $T$ has a
fixed point.

In case B, there exists $n_{0}\in \mathbb{N}$ such that%
\begin{equation}
d(x_{n},T(u))\leq d(u,T(u))\text{,}  \tag{2}
\end{equation}%
for every $n\in \mathbb{N}$, $n\geq n_{0}$. Since $d(u,T(u))\leq
s(d(u,x_{n+1})+d(x_{n+1},T(u)))$, i.e. $\frac{d(u,T(u))}{s}-d(u,x_{n+1})\leq
d(x_{n+1},T(u))$, we get that%
\begin{equation*}
\frac{d(u,T(u))}{s}-d(u,x_{n+1})\leq d(x_{n+1},T(u))\overset{(1)}{\leq }
\end{equation*}%
\begin{equation*}
\leq \alpha \max \{d(x_{n},u),cd(x_{n},x_{n+1}),cd(u,T(u)),\frac{d}{2}%
(d(x_{n},T(u))+d(u,x_{n+1})))\}\overset{(2)}{\leq }
\end{equation*}%
\begin{equation*}
\leq \alpha \max \{d(x_{n},u),cd(x_{n},x_{n+1}),cd(u,T(u)),\frac{d}{2}%
(d(u,T(u))+d(u,x_{n+1})))\}\text{,}
\end{equation*}%
for every $n\in \mathbb{N}$, $n\geq n_{0}$. By passing to limit as $%
n\rightarrow \infty $ in the above inequality, we obtain that%
\begin{equation*}
d(u,T(u))\leq \alpha s\max \{0,0,cd(u,T(u)),\frac{d}{2}d(u,T(u))\}=\alpha
s\max \{c,\frac{d}{2}\}d(u,T(u))\text{.}
\end{equation*}%
As $\alpha s\max \{c,\frac{d}{2}\}<1$ (see ii)), we infer that $d(u,T(u))=0$%
, so $u\in T(u)$, i.e. $T$ has a fixed point. $\square $

\newpage

\textbf{4. Remarks and comments}

\bigskip

\textbf{I}. Let us recall the following result (see Lemma 3.1 from [28]):

\bigskip

\textbf{Lemma 4.1}. \textit{Every sequence} $(x_{n})_{n\in 
%TCIMACRO{\U{2115} }%
%BeginExpansion
\mathbb{N}
%EndExpansion
}$ \textit{of elements from a }$b$\textit{-metric space }$(X,d)$ \textit{is
Cauchy provided that:}

\textit{i) there exists }$\gamma \in \lbrack 0,1)$ \textit{such that }%
\begin{equation*}
d(x_{n+1},x_{n})\leq \gamma d(x_{n},x_{n-1})\text{,}
\end{equation*}%
\textit{\ for every }$n\in \mathbb{%
%TCIMACRO{\U{2115} }%
%BeginExpansion
\mathbb{N}
%EndExpansion
}$\textit{;}

\textit{ii) }$s\gamma <1$\textit{.}

\bigskip

Obviously our Lemma 2.2 is a generalization of the above Lemma which is the
corner stone of the results from [16], [18], [19], [22] and [28].

\bigskip

\textbf{II}. The following definition is inspired by the definition of a
multi-valued weakly Picard operator in the setting of a metric space from
[7].

\bigskip

\textbf{Definition}. \textit{A function }$T:X\rightarrow \mathcal{B}(X)$%
\textit{, where }$(X,d)$\textit{\ is a }$b$\textit{-metric space, is called
a multi-valued weakly Picard operator if for each }$x\in X$\textit{\ and
each }$y\in T(x)$\textit{\ there exists a sequence }$(x_{n})_{n\in \mathbb{N}%
}$\textit{\ such that:}

\textit{i) }$x_{0}=x$\textit{\ and }$x_{1}=y$\textit{;}

\textit{ii) }$x_{n+1}\in T(x_{n})$\textit{\ for every }$n\in \mathbb{N}$%
\textit{;}

\textit{iii) the sequence }$(x_{n})_{n\in \mathbb{N}}$\textit{\ is
convergent and its limit is a fixed point of }$T$\textit{.}

\bigskip

Let us mention that Theorems 3.1., 3.2 and 3.3 provide sufficient conditions
for a function $T$ to be multi-valued weakly Picard operator.

\bigskip

\textbf{III}. For $c=d=0$ in Theorem 3.3 we obtain Theorem 5 from [21], i.e.
Nadler's contraction principle for set-valued functions.

\bigskip

\textbf{IV}. Let us recall the following result (see Theorem 2.2 from [4])
which is a generalization of Theorem 1.2 from [2] which improves Theorem 3.3
from [15], Corollary 3.3. from [26], Corollary 4.3 from [28] and Theorem 1
from [11]:

\bigskip

\textbf{Theorem 4.1}. \textit{A function} $T:X\rightarrow \mathcal{B}(X)$%
\textit{, where }$(X,d)$\textit{\ is a complete }$b$\textit{-metric space,
has a fixed point, provided that it satisfies the following two conditions:}

\textit{i) there exists }$\alpha \in \lbrack 0,1)$\textit{\ such that} 
\begin{equation*}
h(T(x),T(y))\leq \alpha \max
\{d(x,y),d(x,T(x)),d(y,T(y)),d(x,T(y)),d(y,T(x))\}\text{,}
\end{equation*}%
\textit{for all }$x,y\in X$\textit{;}

\textit{ii) }$\alpha \leq \frac{1}{s+s^{2}}$.

\bigskip

Our Theorem 3.3 is a generalization of Theorem 4.1.

Indeed, on the one hand, i) from Theorem 4.1 is a particular case of i) from
Theorem 3.3 since 
\begin{equation*}
N_{1,1}(x,y)=\max \{d(x,y),d(x,T(x)),d(y,T(y)),\frac{1}{2}%
(d(x,T(y))+d(y,T(x)))\}
\end{equation*}%
and%
\begin{equation*}
\frac{1}{2}(d(x,T(y))+d(y,T(x)))\leq \max \{d(x,T(y)),d(y,T(x))\}\text{,}
\end{equation*}%
for all $x,y\in X$. One the other hand, ii) from Theorem 4.1 is a particular
case of ii) from Theorem 3.3 since $\alpha \leq \frac{1}{s+s^{2}}\Rightarrow
\max \{\alpha cs,\alpha ds\}<1$.

\bigskip

Now let us present a situation when Theorem 3.3 is applicable, but Theorem
4.1 is not.

We consider the $b$-metric space $(\mathbb{R},d)$, where $d(x,y)=(x-y)^{2}$
for all $x,y\in \mathbb{R}$, for which $s=2$ and the function $f:\mathbb{R}%
\rightarrow \mathcal{B}(\mathbb{R})$ given by $f(x)=\{\frac{9}{10}x\}$ for
every $x\in \mathbb{R}$. On the one hand, Theorem 3.3. is applicable taking $%
c=d=0$ and $\alpha =\frac{9}{10}$. On the other hand Theorem 4.1 is not
applicable since i) implies $\frac{9}{10}\leq \alpha $ and ii) implies $%
\alpha \leq \frac{1}{6}$.

\bigskip

\textbf{References}

\bigskip

[1] A. Aghajani, M. Abbas and J.R. Roshan, Common fixed point of generalized
weak contractive mappings in partially ordered $b$-metric spaces, Math.
Slovaca, \textbf{64} (2014), 941-960.

[2] A. Amini-Harandi, Fixed point theory for set-valued quasi-contraction
maps in metric spaces, Appl. Math. Lett., \textbf{24} (2011), 1791-1794.

[3] T.V. An, L.Q. Tuyen and N.V. Dung, Stone-type theorem on $b$-metric
spaces and applications, Topology Appl., \textbf{185/186} (2015), 50-64.

[4] H. Aydi, M.F. Bota, E. Karapinar and S. Mitrovi\'{c}, A fixed point
theorem for set-valued quasi-contractions in $b$-metric spaces, Fixed Point
Theory Appl., 2012, 2012:88.

[5] I. A. Bakhtin, The contraction mapping principle in quasimetric spaces,
Funct. Anal., Unianowsk Gos. Ped. Inst., \textbf{30} (1989), 26-37.

[6] V. Berinde, Generalized contractions in quasimetric spaces, Seminar on
Fixed Point Theory, 1993, 3-9.

[7] M. Berinde and V. Berinde, On a general class of multi-valued weakly
Picard mappings, J. Math. Anal. Appl., \textbf{326} (2007), 772-782

[8] M. Boriceanu, M. Bota and A. Petru\c{s}el, Multivalued fractals in $b$%
-metric spaces, Cent. Eur. J. Math., \textbf{8} (2010), 367-377.

[9] M. Boriceanu, A. Petru\c{s}el and A.I. Rus, Fixed point theorems for
some multivalued generalized contraction in $b$-metric spaces, Int. J. Math.
Stat., \textbf{6} (2010), 65-76.

[10] M. Bota, A. Moln\'{a}r and C. Varga, On Ekeland's variational principle
in $b$-metric spaces, Fixed Point Theory, \textbf{12} (2011), 21-28.

[11] L. Ciri\'{c}, A generalization of Banach's contraction principle, Proc.
Amer. Math. Soc., \textbf{45} (1974), 267-273.

[12] C. Chifu and G. Petru\c{s}el, Fixed points for multivalued contractions
in $b$-metric spaces with applications to fractals, Taiwanese J. Math., 
\textbf{18} (2014), 1365-1375.

[13] S. Czerwik, Contraction mappings in $b$-metric spaces, Acta Math.
Inform. Univ. Ostraviensis, \textbf{1} (1993), 5-11.

[14] S. Czerwik, Nonlinear set-valued contraction mappings in $b$-metric
spaces, Atti Sem. Mat. Fis. Univ. Modena, \textbf{46} (1998), 263-276.

[15] P.Z. Daffer and H. Kaneko, Fixed points of generalized contractive
multi-valued mappings, J. Math. Anal. Appl., \textbf{192} (1995), 655-666.

[16] A.K. Dubey, R. Shukla and R.P. Dubey, Some fixed point results in $b$%
-metric spaces, Asian Journal of Mathematics and Applications, 2014, Article
ID ama0147.

[17] M.A. Khamsi and N. Hussain, KKM mappings in metric type spaces,
Nonlinear Anal., \textbf{73} (2010), 3123-3129.

[18] M. Kir, H. Kizitune, On some well known fixed point theorems in $b$%
-metric spaces, Turkish Journal of Analysis and Number Theory, \textbf{1}
(2013), 13-16.

[19] P. K. Mishra, S. Sachdeva and S.K. Banerjee, Some fixed point theorems
in $b$-metric space, Turkish Journal of Analysis and Number Theory, \textbf{2%
} (2014), 19-22.

[20] S. K. Mohanta, Some fixed point theorems using $wt$-distance in $b$%
-metric spaces, Fasc. Math., \textbf{54} (2015), 125-140.

[21] S.B. Nadler, Multi-valued contraction mappings, Pac. J. Math., \textbf{%
30} (1969), 475-488.

[22] H.N. Nashine and Z. Kadelburg, Cyclic generalized $\varphi $%
-contractions in $b$-metric spaces and an application to integral equations,
Filomat, \textbf{28} (2014), 2047-2057.

[23] M.O. Olatinwo, A fixed point theorem for multi-valued weakly Picard
operators in $b$-metric spaces, Demonstratio Math., \textbf{42} (2009),
599-606.

[24] M. P\u{a}curar, Sequences of almost contractions and fixed points in $b$%
-metric spaces, An. Univ. Vest Timi\c{s}. Ser. Mat.-Inform., \textbf{48}
(2010), 125-137.

[25] J.R. Roshan, N. Hussain, S. Sedghi and N. Shobkolaei, Suzuki-type fixed
point results in $b$-metric spaces, Math. Sci. (Springer), \textbf{9}
(2015), 153--160.

[26] B.D. Rouhani and S. Moradi, Common fixed point of multivalued
generalized $\varphi $-weak contractive mappings, Fixed Point Theory Appl.,
2010, Article ID 708984.

[27] M. Sarwar and M.U. Rahman, Fixed point theorems for Ciric's and
generalized contractions in $b$-metric spaces, International Journal of
Analysis and Applications, \textbf{7 }(2015), 70-78.

[28] S.L. Singh, S. Czerwik, K. Kr\'{o}l and A. Singh, Coincidences and
fixed points of hybrid contractions, Tamsui Oxf. J. Math. Sci., \textbf{24}
(2008), 401-416.

[29] S. Shukla, Partial $b$-metric spaces and fixed point theorems,
Mediterr. J. Math., \textbf{11} (2014), 703-711.

[30] H. Yingtaweesittikul, Suzuki type fixed point for generalized
multi-valued mappings in $b$-metric spaces, Fixed Point Theory and
Applications, 2013, 2013:215.

\bigskip

University of Bucharest

Faculty of Mathematics and Computer Science

Str. Academiei\ 14, 010014 Bucharest, Romania

E-mail: miculesc@yahoo.com, mihail\_alex@yahoo.com

\end{document}